\documentclass{amsart}
\input epsf

\usepackage{amsfonts,amsthm,amsmath,amssymb,latexsym}
\usepackage{graphicx,color}
\usepackage[all]{xy}

\def\L{\mathcal{L}}
\def\D{{\bf D}}

\def\G{\mathcal{G}}

\def\S{\mathcal{S}}

\begin{document}

\author{Dragomir \v Sari\' c}

\address{Institute for Mathematical Sciences, Stony Brook University,
Stony Brook, NY 11794-3660} \email{saric@math.sunysb.edu}

\title[Earthquakes and Thurston's boundary for the solenoid]{Earthquakes
and Thurston's boundary for the Teichm\"uller space of
the universal hyperbolic solenoid}


\date{\today}

\begin{abstract}
A measured laminations on the universal hyperbolic solenoid $\S$ is,
by our definition, a leafwise measured lamination with appropriate
continuity for the transverse variations. An earthquakes on the
universal hyperbolic solenoid $\S$ is uniquely determined by a
measured lamination on $\S$; it is a leafwise earthquake with the
leafwise earthquake measure equal to the leafwise measured
lamination. Leafwise earthquakes fit together to produce a new
hyperbolic metric on $\S$ which is transversely continuous and we
show that any two hyperbolic metrics on $\S$ are connected by an
earthquake. We also establish the space of projective measured
lamination $PML(\S )$ as a natural Thurston-type boundary to the
Teichm\"uller space $T(\S )$ of the universal hyperbolic solenoid
$\S$. The (baseleaf preserving) mapping class group $MCG_{BLP}(\S )$
acts continuously on the closure $T(\S )\cup PML(\S )$ of $T(\S )$.
Moreover, the set of transversely locally constant measured
laminations on $\S$ is dense in $ML(\S )$.
\end{abstract}

\maketitle

\section{Introduction}

The universal hyperbolic solenoid $\S$ is the inverse limit of the
system of all finite sheeted unbranched pointed covers of a compact
surface of genus greater than $1$. Sullivan \cite{Sul} introduced
the solenoid $\S$ as the ``universal compact surface'', i.e. the
universal object in the category of finite unbranched covers of a
compact surface. It turns out that $\S$ has a rich deformation
theory, i.e. the Teichm\"uller space $T(\S )$ is highly non-trivial.
In fact, $T(\S )$ is a first example of a Teichm\"uller space which
is separable but not finite dimensional. (Recall that Teichm\"uller
spaces of compact Riemann surfaces with possibly finitely many
points removed are finite-dimensional complex manifolds, while
Teichm\"uller spaces of geometrically infinite Riemann surfaces are
non-separable infinite-dimensional complex Banach manifolds.)

\vskip .2 cm

Sullivan \cite{Sul} started the study of the complex structure and
the Teichm\"uller metric on the Teichm\"uller space $T(\S )$ of the
universal hyperbolic solenoid $\S$. The universal hyperbolic
solenoid $\S$ has a transverse measure (unlike most of laminations)
which is utilized in \cite{Sa4} to continue the investigation of the
Teichm\"uller metric on $T(\S )$. The results in \cite{Sa4} are used
in \cite{EMS} to show that generic points in $T(\S )$ do not have
Teichm\"uller-type Beltrami coefficient representatives which
sharply contrasts the situation for finite surfaces (all points have
Teichm\"uller type representatives) and for infinite surfaces (a
dense, open subset of the Teichm\"uller space does). Thus, most
points in $T(\S )$ are not connected to the basepoint by a nice
Teichm\"uller geodesic obtained by stretching horizontal and
shrinking vertical foliation of a holomorphic quadratic differential
on $\S$ unlike for Riemann surfaces (see \cite{Lak} for the
statement for infinite surfaces).

\vskip .2 cm

By results of Candel \cite{Ca}, a conformal structure on $\S$
contains a unique hyperbolic metric. This paper starts the
investigation of the Teichm\"uller space $T(\S )$ of the solenoid
$\S$ using hyperbolic structures on $\S$. (An approach to studying
the decorated Teichm\"uller space of the related punctured solenoid
via hyperbolic structures is made in \cite{PS}, \cite{BPS}. This
approach is specific to the punctured solenoid while the results in
the paper hold for both the universal hyperbolic solenoid and the
punctured solenoid.)

\vskip .2 cm

 Our first result concerns the notion of an earthquake between two
 hyperbolic solenoids. An earthquake on a hyperbolic surface is a
 piecewise isometric bijective (not necessarily continuous)
 map from the hyperbolic surface to another hyperbolic surface. The support of an
 earthquake is a geodesic lamination (called the earthquake fault \cite{Th}) along which the quaking
 (i.e. the discontinuity) appears. The restriction to each stratum (i.e. a leaf of the support or a connected
 component of the complement of the support) of an earthquake is
 an isometry, and each stratum is moved to the left when viewed
 from any other stratum. Given an earthquake on a hyperbolic surface, there exists
 a unique transverse measure (called an earthquake measure) to its support \cite{Th} which
 determines the earthquake. (An earthquake measure is identified with a positive
 Radon measure, called a measured lamination, on the space of geodesics of a surface whose support
 is the support geodesic lamination and the measure of a bunch of geodesics in the support
 is given by the measure of an arc which intersects them and which does
 not intersect
 other geodesics in the support.)
 The main result concerning the
 earthquakes on hyperbolic surfaces is that any two hyperbolic metrics are
 connected by an earthquake \cite{Th}.

 \vskip .2 cm

 We introduce a proper notion of an earthquake measure on the
 universal hyperbolic solenoid $\S$. An earthquake measure on $\S$
 is an assignment of an earthquake measure on each leaf (isometric to a hyperbolic plane) of
 $\S$ such that the measures vary continuously for the transverse variations in an appropriate
 topology (see Definition 4.2 for more details). The support of
 leafwise earthquake measures are leafwise geodesic laminations
 which do not vary continuously for the transverse variations (see
 Example 4.3). However, the continuity of measures (in the
 appropriate Fr\'echet topology) guarantees that leafwise
 earthquakes on $\S$ piece together a new hyperbolic structure on
 $\S$ which is continuous for the transverse variations.
 Therefore, earthquake measures on $\S$ produce an
 earthquake of $\S$. We establish the transitivity statement for
 earthquakes on hyperbolic structures of the universal hyperbolic
 solenoid $\S$ analogous to the case of hyperbolic surfaces.

 \vskip .2 cm

 \paragraph{\bf Theorem 5.1} {\it A measured lamination $\mu$ on a solenoid $X$ with
 an arbitrary hyperbolic metric gives an earthquake map $E_{\mu}$ of $X$ into another
 solenoid $Y$ with a hyperbolic metric such that there exists a (differentiable)
 quasiconformal map $f:X\to Y$ whose extension to the boundary of leaves coincides
 with the extension of $E_{\mu}$. Any two points in the Teichm\"uller space $T(\S )$
 of the universal hyperbolic
 solenoid $\S$ are connected by a unique earthquake along a measured
 lamination.}

 \vskip .2 cm

To prove Theorem 5.1, we showed that if two quasisymmetric maps of
the unit circle $S^1$ are close then their corresponding earthquake
measures are close (in the Fr\'echet topology). The opposite is
false by an easy example. Therefore, the earthquake map from the
space of bounded measured laminations in the unit disk onto the
Teichm\"uller space is bijective but not a homeomorphism for the
Fr\'echet topology on measured lamination. In the case of compact
surfaces, the earthquake map is a homeomorphism (see \cite{Ker})
when measured lamination are equipped with weak* topology (which is
equivalent to the Fr\'echet topology for compact surfaces). We
establish similar result for the solenoid $\S$.

\vskip .2 cm

\paragraph{\bf Corollary 5.2} {\it The earthquake map which assigns
to each bounded measured lamination on the universal hyperbolic
solenoid $\S$ the corresponding marked hyperbolic solenoid is a
homeomorphism between the space $ML(\S )$ of bounded measured
laminations and the Teichm\"uller space $T(\S )$.}

\vskip .2 cm

Thurston \cite{Thu2},\cite{FLP} introduced a natural boundary to the
Teichm\"uller space of a compact surface by ``adding at infinity''
the space of projective measured laminations. The mapping class
group acts continuously on the closure and there is a classification
of its elements according to their action on the boundary. Bonahon
\cite{Bo1} gave an alternative description of the Thurston's
boundary to the Teichm\"uller space of a compact surface using the
Liouville map which embeds the Teichm\"uller space into the space of
measures on the space of geodesics of the surface. The boundary
points at infinity are asymptotic rays to the image of the
Teichm\"uller space. We used (see \cite{Sa3}) the idea of the
Liouville embedding to give a Thurston-type boundary to the
Teichm\"uller space of any (possibly geometrically infinite) Riemann
surface. We extend this result to the Teichm\"uller space $T(\S )$
of the universal hyperbolic solenoid $\S$.

\vskip .2 cm

Biswas, Nag and Mitra \cite{BNM} introduced the direct limit of projective
measured laminations on the compact surfaces as boundary at infinity of
the direct limit of
Teichm\"uller spaces of compact surfaces covering a fixed compact surface of
genus at least two. Since $T(\S )$ contains
as a dense subset the above direct limit
of Teichm\"uller spaces of compact surfaces, they remarked that the Thurston's boundary
for $T(\S )$ should be a completion of the union of the projective measured laminations on all
compact surfaces. We give a proper analytical description of the Thurston's boundary
answering their question about the completion.
 The main point is to properly
define the continuity for the transverse variations of various
spaces of measures and distributions on the space $\G (\S )$ of
geodesics on the universal hyperbolic solenoid $\S$. We establish
this goal using the Fr\'echet topology on the (``enveloping'') space
of H\"older distributions $H(\S )$ (see Section 6).

\vskip .2 cm

\paragraph{\bf Theorem 6.2} {\it The Liouville map
$\mathcal{L}_{\S}:T(\S )\to H(\S )$ is a homeomorphisms onto
its image. The set of asymptotic rays to $\mathcal{L}_{\S}(T(\S ))$ is
homeomorphic to the space of projective measured laminations on $\S$. The
baseleaf preserving mapping class group $MCG_{BLP}(\S )$ acts continuously
on the closure $T(\S )\cup PML(\S )$ of the Teichm\"uller space $T(\S )$
of the universal hyperbolic solenoid $\S$.}

\vskip .2 cm

A lift of a measured lamination on a compact surface of genus at least two to the universal
hyperbolic solenoid $\S$ is a measured lamination on $\S$. Such measured lamination is locally
constant for the transverse variations and is called a transversely locally constant (TLC)
measured lamination. We showed that
an arbitrary measured lamination on $\S$ is the limit in the Fr\'echet topology of the TLC
measured laminations.

\vskip .2 cm

\paragraph{\bf Theorem 6.3} {\it The subset of all measured
lamination on the universal hyperbolic solenoid $\S$ which are locally transversely constant is dense
in the space of all measured laminations $ML(\S )$ on $\S$ for the Fr\'echet topology.}

\vskip .2 cm

In Section 7, we introduce the space of compactly supported measured lamination $PML_0(\S_p )$
on the punctured solenoid $\S_p$
and extend the transitivity of earthquakes and Thurston's boundary for the punctured solenoid
by replacing $ML(\S )$ with $ML_0(\S_p)$.

\section{Preliminaries}

We recall several definitions: the universal hyperbolic solenoid, earthquakes of the unit disk,
the Fr\'echet topology on the measures on the space of geodesics of the unit disk.

\subsection{The universal hyperbolic solenoid}

Let $(S,x)$ be a fixed compact surface of genus at least two with
the basepoint $x\in S$. Consider all finite sheeted unbranched
covers $(S_i,x_i)$ by compact surfaces with basepoints such that
the covering maps $\pi_i:S_i\to S$ satisfy $\pi_i(x_i)=x$. There
is a natural partial ordering $\leq$ on the set of all such
coverings. Namely, $(S_i,x_i)\leq (S_j,x_j)$ if there exists a
finite sheeted unbranched covering map $\pi_{j,i}:S_j\to S_i$,
$\pi_{j,i}(x_j)=x_i$, such that $\pi_i\circ \pi_{j,i} =\pi_j$. The
set of all coverings is inverse directed, i.e. given two coverings
of $S$ there exists a third covering of $S$ which is larger than
the two (see \cite{Sul}, \cite{Od}, \cite{MS}, \cite{Sa4}).
Sullivan \cite{Sul} introduced the universal hyperbolic solenoid
$\S$ as follows.

\vskip .2 cm

\paragraph{\bf Definition 2.1} The {\it universal hyperbolic
solenoid} $\S$ is the inverse limit (for the above partial ordering)
of the directed system of all finite sheeted unbranched
covers of a fixed compact surface of genus at least two.

\vskip .2 cm

The inverse limit $\S$ is independent of the base surface (i.e. two inverse limits
with two base surfaces of genus at least two
are homeomorphic). Thus it is called the universal hyperbolic solenoid $\S$.

\vskip .2 cm

We give an equivalent definition of the universal hyperbolic
solenoid $\S$ \cite{Od}. Let $G$ be a Fuchsian group uniformizing
$S$, i.e. $S$ is homeomorphic to $\D /G$, where $\D$ is the unit
disk. Let $G_n$ be the intersection of all subgroups of $G$ with
index at most $n$. Then $G_n$ is of finite index in $G$. We define
{\it profinite metric} on $G$ (see \cite{Od}) by
$$
dist(A,B)=e^{-n},
$$
for $A,B\in G$, where $AB^{-1}\in G_n$ and $AB^{-1}\notin G_{n+1}$. The completion of
$G$ in the profinite metric $dist$ is called the {\it profinite group completion } $\hat{G}$.
The completion $\hat{G}$ is a compact, topological group homeomorphic to Cantor set. The group $G$
is a dense subgroup of $\hat{G}$. We define an action of $G$ on the product $\D\times\hat{G}$ by
$$
A(z,t):=(A(z),tA^{-1}),
$$
where $z\in\D$ and $t\in\hat{G}$, the action of $A$ on $\D$ is
just a M\"obius map and $A$ acts on $\hat{G}$ because $G$ lies
inside $\hat{G}$. The universal hyperbolic solenoid $\S$ is
homeomorphic to the quotient $(\D\times\hat{G})/G$ (see
\cite{Od}). From this description of $\S$, it is easy to see that
$\S$ is a compact space which is locally homeomorphic to a 2-disk
times Cantor set. Moreover, $\S$ fibers over $\D /G$ with fibers
Cantor sets isomorphic to $\hat{G}$, and the restriction of the fiber map to each leaf is
the universal covering of $S\equiv \D /G$. The same is true for
any finite cover of $S$ \cite{Odd} by replacing $G$ with the covering group.
Each $(\D\times\{ t\}
)/G\subset\S$, for $t\in\hat{G}$, is a path component, called a
{\it leaf} of $\S$. Each leaf is homeomorphic to the unit disk and
it is dense in $\S$. Thus, each leaf
is simply connected but the restriction of the topology on $\S$ to
a leaf is not the standard topology on the unit disk. We define the
universal cover of $\S$ by ``straightening'' the topology on
leaves.

\vskip .2 cm

\paragraph{\bf Definition 2.2} The {\it universal cover} for the
universal hyperbolic solenoid $\S$ is given by
$$
\pi :\D\times\hat{G}\to (\D\times\hat{G})/G,
$$
where $\pi$ is the quotient map for the action of $G$.

\vskip .2 cm

A {\it complex structure} on the universal hyperbolic solenoid $\S$
is a collection of local charts such that the transition maps are
holomorphic in the disk direction and they vary continuously (for
topology of uniform convergence) in the transverse (Cantor)
direction of the local charts \cite{Sul}. Complex structures on $\S$
are in one to one correspondence with conformal structures (which
are continuous for the variations in the transverse direction) on
$\S$ by the continuous dependence on the parameters of the solution
of Beltrami equation. Any conformal structure on $\S$ contains a
unique hyperbolic metric which is continuous for the variations in
the transverse direction (see \cite{Ca}). A hyperbolic metric (or a
complex structure) on $\S$ which is transversely locally constant
for some choice of charts on $\S$ is called a TLC hyperbolic metric
(or a TLC complex structure) on $\S$ \cite{Sul}. Any TLC hyperbolic
metric (or complex structure) on $\S$ is obtained by lifting a
hyperbolic metric (or a complex structure) from a finite cover of
$S$ to $\S$ \cite{NS}. Note that by identifying $\S$ with
$(\D\times\hat{G})/G$ we fix a TLC complex structure on $\S$ coming
from Riemann surface $\D /G$. A {\it (differentiable) quasiconformal
map} $f:\S\to X$ from the fixed TLC complex solenoid $\S$ to an
arbitrary complex solenoid $X$ is a homeomorphism which is
$C^{\infty}$-differentiable in the disk direction in local charts
and varies continuously in the transverse direction for
$C^{\infty}$-topology on $C^{\infty}$-maps \cite{Sul}, \cite{Sa4}.
We note that the use of differentiable quasiconformal maps as
opposed to only quasiconformal maps is necessary in order for
compositions of quasiconformal maps to be continuous in the
transverse direction. We define the Teichm\"uller space $T(\S )$ of
the universal hyperbolic solenoid $\S$.

\vskip .2 cm

\paragraph{\bf Definition 2.3} The {\it Teichm\"uller space} $T(\S )$ of
the universal hyperbolic solenoid $\S$ consists of all
quasiconformal maps $f:\S\to X$ up to an equivalence. Two
quasiconformal maps $f:\S\to X$ and $g:\S\to X_1$ are equivalent if
there exists a conformal map $c:X\to X_1$ such that $g^{-1}\circ
c\circ f:\S\to\S$ is homotopic to the identity. The equivalence
class of the identity $id:\S\to\S$ is called the {\it basepoint} of
$T(\S )$.

\subsection{Earthquakes in the unit disk}

We define earthquakes of the unit disk $\D$ and recall their main
properties. A {\it geodesic lamination} in the unit disk $\D$ is a
closed subset of $\D$ which is foliated by geodesics for the
hyperbolic metric on $\D$, or equivalently, it is a closed subset of
the space $\mathcal{G}(\D )$ of geodesics in $\D$ such that no two
geodesics in the subset intersect in $\D$. Some familiar examples of
geodesic laminations in $\D$ are: a set of finitely many
non-intersecting geodesics in $\D$; a countable, discrete set of
non-intersecting geodesics; a foliation of $\D$ by geodesics.

\vskip .2 cm

\paragraph{\bf Definition 2.4} An {\it earthquake measure} on the
unit disk is a positive Radon measure on the space of geodesic $\G
(\D )$ of the unit disk $\D$ whose support is a geodesic lamination.

\vskip .2 cm

Note that $\G (\D )$ is homeomorphic to $(S^1\times
S^1-diag)/\mathbb{Z}_2$ by mapping a geodesic in $\D$ to the
unordered pair of its ideal endpoints in $S^1$. Thurston \cite{Th}
introduced earthquakes as follows.

\vskip .2 cm

\paragraph{\bf Definition 2.5} An {\it earthquake} $E:\D\to\D$ of the unit
disk $\D$ is a bijective map which maps a fixed geodesic lamination $\lambda$ (called the {\it support} of $E$)
in $\D$  onto another geodesic lamination $\lambda'$.
A geodesic from $\lambda$ or a connected component of $\D -\lambda$ is called a {\it stratum}
of $E$. The restriction of the earthquake $E$ to each stratum is a hyperbolic isometry with
the additional property that for any two strata $A,B$ of $E$, the {\it comparison isometry}
$$
E|_B\circ (E|_A)^{-1}
$$
is a hyperbolic translation whose axis separates $A$ from $B$, and which translates $B$ to the left
as seen from $A$.

\vskip .2 cm

Each earthquake $E:\D\to\D$ continuously extends to a
homeomorphism of $\partial\D\equiv S^1$, which we denote by
$E|_{S^1}:S^1\to S^1$ \cite{Th}. An important theorem due to Thurston
is that each orientation
preserving homeomorphism of $S^1$ is obtained as the extension to
$S^1$ of an earthquake \cite{Th}. Given an earthquake $E$ of $\D$,
there exists a unique corresponding earthquake measure $\mu$ on
$\G (\D )$ supported on $\lambda$ determined by the following
condition. Consider a subset of $\lambda$ consisting of geodesics
which intersect a closed arc $I$. Choose finitely many strata of
$E$ intersecting $I$. The $\mu$ measure of the subset is
approximated by the sum of translation lengths between comparison
isometries of adjacent strata (of the above chosen finitely many
strata of $E$ intersecting the arc $I$) when the distance between
adjacent strata goes to zero \cite{Th}. An alternative description
is to consider the measure $\mu$ to be a family of measures on arcs $I$ in
$\D$ which are invariant under homotopies of arcs preserving the leaves of
$\lambda$. If an earthquake measure $\mu$ corresponds to an
earthquake as above, we denote the corresponding earthquake by
$E_{\mu}$. Two homeomorphisms have the same corresponding
earthquake measures if and only if they differ by a
post-composition with a hyperbolic isometry of $\D$ \cite{Th}.

\vskip .2 cm

We say that an earthquake measure $\mu$ is {\it bounded} if the {\it norm}
$\|\mu\|$ satisfies
$$
\|\mu\| :=\sup_I\mu (I)<\infty ,
$$
where the supremum is over all geodesic arcs in $\D$ of length
$1$. If $\mu$ is a bounded earthquake measure, then there exists
earthquake $E_{\mu}$ corresponding to $\mu$. A homeomorphism
$h:S^1\to S^1$ is quasisymmetric if and only if the corresponding
earthquake measure $\mu$ is bounded, where $E_{\mu}|_{S^1}=h$
\cite{Sa}, \cite{Sa1}.

\subsection{The Fr\'echet topology}

We recall the definition of the Fr\'echet topology on the space of
H\"older distributions $H(\D )$ of the unit disk $\D$ from
\cite{Sa2}. The space of bounded positive measures on $\G (\D )$,
and, in particular, the space of bounded earthquake measures on $\G
(\D )$ are subsets of $H(\D )$.

\vskip .2 cm

The {\it Liouvile measure} $L$ on $\mathcal{G} (\D)$ is given by
$$
L(K):=\iint_K\frac{d\alpha d\beta}{|e^{i\alpha}-e^{i\beta}|^2}\ ,
$$
where $e^{i\alpha},e^{i\beta}\in S^1$. A {\it box of geodesics} $Q$
is the set of all geodesics in $(a,b)\times (c,d)\subset S^1\times
S^1-diag$, where $a,b,c,d\in S^1$ are different points given in the
counter-clockwise order on $S^1$. Then
$$
L(Q)=\log\Big{|}\frac{(a-c)(b-d)}{(a-d)(b-c)}\Big{|},
$$
and this formula can be used as an alternative definition of Liouville measure.

\vskip .2 cm

We fix $0<\nu\leq 1$.
The space of $\nu$-{\it test functions} $test(\nu )$ consists of all
$\nu$-H\"older continuous functions $(\varphi ,Q)$, $\varphi :\mathcal{G}(\D )\to\mathbb{R}$,
whose support is in a box of geodesics $Q$ with $L(Q)=\log 2$ such
that $\|\varphi\circ\Theta_Q\|_{\nu}\leq 1$, where
$\Theta_Q:(-1,-i)\times (1,i)\mapsto Q$ is a hyperbolic isometry
and $\|\varphi\|_{\nu} :=\max\{\sup_{\mathcal{G}(\D )}|\varphi |,
\sup_{(x,y)\neq (x_1,y_1)}\frac{|\varphi (x,y)-\varphi
(x_1,y_1)|}{d((x,y),(x_1,y_1))^{\nu}}\}$, with $d$ being the angle
metric on $S^1\times S^1$ (see \cite{Sa2},\cite{Sa3}).

\vskip .2 cm

The space of {\it H\"older distribution} $H(\D )$ (see
\cite{Sa3}) of the unit diks $\D$ consists of all linear functionals
$\Psi$ on the space of H\"older continuous functions $\varphi
:\mathcal{G}(\D )\to\mathbb{R}$ with compact support such that
$$
\|\Psi\|_{\nu}:=\sup_{\varphi\in test(\nu )}|\Psi (\varphi )|<\infty
$$
for all $0<\nu\leq 1$. The Fr\'echet topology on $H(\D )$ is
defined using the family of $\nu$-norms above. The topological vector
space $H(\D )$ is metrizable and one metric which gives the Fr\'echet
topology is
$$
dist(\Psi ,\Psi_1 ):=\sum_{n=1}^{\infty}\frac{1}{n^2}\|\Psi -\Psi_1\|_{1/n}.
$$

\section{The convergence of measures in the unit disk}

Denote by $\mathcal{G}_z$, for $z\in\D$, the set of geodesics in
$\D$ which contain $z$. If $z,w\in\D$ then denote by $[z,w]$ the
geodesic arc in $\D$ between $z$ and $w$. If $K$ is a subset of
$\D$, denote by $\mathcal{G}_K$ the set of geodesics of $\D$ which
intersect $K$.

\vskip .2 cm

We showed in \cite{Sa1} that a sequence of homeomorphisms of $S^1$
pointwise converges to a homeomorphism of $S^1$ if and only if the
sequence of earthquake measures, corresponding to the sequence of
homeomorphisms, converges to the earthquake measure of the limit.
More precisely,

\vskip .2 cm

\paragraph{\bf Proposition 3.1} \cite{Sa1} {\it Let $\mu ,\mu_i$ be
uniformly bounded earthquake measures on $\D$, i.e. $\|\mu\|
,\|\mu_i\|\leq M<\infty$. Then $\mu_i\to\mu$ in the weak* topology
as $i\to\infty$ if and only if there exist normalizations of
earthquake maps $E_{\mu_i}|_{S^1},E_{\mu}|_{S^1}$ such that
$E_{\mu_i}|_{S^1}(x)\to E_{\mu}|_{S^1}(x)$ for each $x\in S^1$, as
$i\to\infty$. (Note that $E_{\mu_i}|_{S^1},E_{\mu}|_{S^1}$ are
well-defined up to the post-compositions by isometries and different
normalizations correspond to different choices of isometries.)}

\vskip .2 cm

We consider a sequence of quasisymmetric maps converging to a
quasisymmetric map in the quasisymmetric topology and show that
the corresponding sequence of earthquake measures converges in the Fr\'echet
topology.

\vskip .2 cm

\paragraph{\bf Proposition 3.2} {\it Let $h_n=E_{\mu_n}|_{S^1}$ and
$h=E_{\mu}|_{S^1}$ be quasisymmetric maps such that $h_n\to h$ as
$n\to\infty$, in the quasisymmetric topology. Then
$\|\mu_n-\mu\|_{\nu}\to 0$ as $n\to\infty$, for each $0<\nu\leq
1$.}

\vskip .2 cm

\paragraph{\bf Proof} Assume on the contrary that
$\|\mu_n-\mu\|_{\nu}\geq m>0$, for a fixed $0<\nu\leq 1$. Thus,
there exists $(\varphi_n,Q_n)\in test(\nu )$ such that
$|\mu_n(\varphi_n)-\mu(\varphi )|\geq m>0$, where $L(Q_n)=\log 2$.
Without loss of generality, we assume that $h_n,h$ fix $1,i,-1$.
Let $Q_n':=h_n(Q_n)$ and $Q_n'':=h(Q_n)$. There exist unique
hyperbolic isometries $A_n,A_n',A_n''$ of the unit disk $\D$ such
that $A_n:Q_n\mapsto (1,i)\times (-1,-i)$, $A_n':Q_n'\mapsto
(1,i)\times (-1,q_n')$ and $A_n'':Q_n''\mapsto (1,i)\times
(-1,q_n'')$, for unique $q_n',q_n''\in S^1$.

\vskip .2 cm

Define $\bar{h}_n:=A_n'\circ h_n\circ A_n^{-1}$ and
$\bar{f}_n:=A_n''\circ h\circ A_n^{-1}$. Note that
$\bar{h}_n:1,i,-1,-i\mapsto 1,i,-1,q_n'$ and
$\bar{f}_n:1,i,-1,-i\mapsto 1,i,-1,q_n''$. The normalization of
$\bar{h}_n$ and $\bar{f}_n$ implies that $\bar{h}_n\to \bar{h}$ and
$\bar{f}_n\to \bar{f}$ pointwise on $S^1$, where $\bar{h},\bar{f}$
are quasisymmetric maps as well. (This convergence is a consequence
of pointwise convergence of a family of $K$-quasiconformal maps
normalized to fix three points in $\hat{\mathbb{C}}$. Note that we
can choose quasiconformal extensions of $\bar{h}_n,\bar{f}_n$ to
have the same quasiconformal constant by using barycentric extension
\cite{DE} in the interior and the exterior of the unit circle
$S^1$.)

\vskip .2 cm

Consequently, $\bar{h}_n\circ \bar{f}_n^{-1}\to
\bar{h}\circ\bar{f}^{-1}$ pointwise, as $n\to\infty$. Let $Q$ be
an arbitrary box with $L(Q)=\log 2$. Then $|L(A_n'\circ h_n\circ
h^{-1}\circ (A_n'')^{-1}(Q))-L(Q)|\to 0$ as $n\to\infty$ because
$h_n\to h$ in the quasisymmetric topology and by the invariance of
Liouville measure under hyperbolic isometries. Thus
$\bar{h}\circ\bar{f}^{-1}$ preserves Liouville measure and fixes
$1,i,-1$. Therefore $\bar{h}=\bar{f}$.

\vskip .2 cm

Let $\mu_n':=A_n^{*}(\mu_n)$ and let $\sigma_n:=A_n^{*}(\mu )$.
Then there exists a sequence $(\varphi_n',(1,i)\times (-1,-i))\in
test(\nu )$ such that
\begin{equation}
\label{difference}
|\mu_n'
(\varphi_n')-\sigma_n'(\varphi_n')|\geq m>0.
\end{equation}
Since $\mu_n'$ is the push forward by a hyperbolic isometry of
$\mu_n$ then $\|\mu_n'\| =\|\mu_n\|$. Moreover, since $\mu_n$ are
earthquake measures for $h_n$ and $h_n$ converges in the
quasisymmetric topology, it follows that $\mu_n$ are uniformly
bounded measures (and the same holds for $\mu_n'$). The sequence
$\sigma_n$ is also uniformly bounded because it is the push
forward of a single measure by hyperbolic isometries. Both
sequences $\mu_n'$ and $\sigma_n$ converge to bounded earthquake
measures $\mu'$ and $\sigma$ such that $\bar{h}=E_{\mu'}|_{S^1}$
and $\bar{f}=E_{\sigma}|_{S^1}$ by Proposition 3.1. By
(\ref{difference}), we conclude that $\mu'\neq\sigma$. (To see
that $\mu'\neq\sigma$, note that the weak* convergence on a fixed
box $(1,i)\times (-1,-i)$ is equivalent to the uniform convergence
with respect to all $(\varphi ,(1,i)\times (-1,-i))\in test(\nu
)$. See \cite{Sa2} for details.) But this is a contradiction with
$\bar{h}=\bar{f}$ by the uniqueness of earthquake measures
\cite{Th}. $\Box$

\vskip .2 cm

We remark that the converse of Proposition 3.2 is not true. This is
easily seen by an example. Take a fixed geodesic with a positive
weight as one lamination. Take a convergent sequence in Fr\'echet
topology to consists of geodesics sharing exactly one endpoint with
the above geodesic and take the same positive weight. It is obvious
that the extension of the earthquakes to $S^1$ corresponding to the
sequence does not converge to the extension of the earthquake to
$S^1$ corresponding to the limit in the quasisymmetric topology.
Note that they do converge pointwise. This is in contrast with the
statement in Proposition 3.1 which gives the equivalence. However,
if we restrict ourselves to the earthquakes on compact surfaces then
the equivalence holds. One of the main results in the next two
sections is that the equivalence holds for the universal hyperbolic
solenoid as well.

\section{Measured laminations on the universal hyperbolic
solenoid}

Recall that a leaf of the universal hyperbolic solenoid $\S$
intersects any local chart countably many times. Each intersection
is a 2-disk, which is called a {\it local leaf}. Given two local leaves
of two global leaves, there exists an identification isometry of
the global leaves given as follows. Since the hyperbolic metric in
the local charts is continuous for the trivial (vertical)
identification, we can choose two points sitting one above the
other and two unit tangent vectors based at the points whose
directions get vertically identified. (The two vectors are not
necessarily vertically identified because the hyperbolic metrics
are not constant in the transverse direction.) The isometric
identification is uniquely determined by requiring to map one point onto
the other other such
that the unit vector is mapped onto the unit vector. The
identification depends on the chart and the choice of two points
while the choice of tangent vectors does not affect it.

\vskip .2 cm

Fix one local leaf $l$ and consider a sequence of local leaves $l_n$
approaching $l$. Suppose we choose two different isometric
identifications $f_n:\tilde{l}\to \tilde{l}_n$ and $g_n:\tilde{l}\to
\tilde{l}_n$ of the global leaves $\tilde{l},\tilde{l}_n$ containing
local leaves $l,l_n$ (the identifications differ by the choice of
points in the local leaves). Then $g_n^{-1}\circ f_n$ is an isometry
of $\tilde{l}$ which converges to the identity as $n\to\infty$
because of the continuity in the transverse direction of the
hyperbolic metrics. This implies that any two identifications of two
global leaves differ by an isometry which is close to the identity
when corresponding local leaves are close. Therefore, it makes sense
to compare objects (preserved by isometries) on two nearby leaves as
well as maps from leaves.

\vskip .2 cm

\paragraph{\bf Definition 4.1} A (transversely continuous)
{\it geodesic lamination} on the universal hyperbolic solenoid $\S$
is an assignment of a geodesic lamination to each leaf which is
continuous (for Hausdorff distance between closed subset of
$\mathcal{G}(\D )$ defined using the angle metric $d$ on $S^1\times
S^1$) with respect to the transverse variations given by each local
chart as above.

\vskip .1 cm

Namely, for any local chart, we consider the isometric
identifications as above. The geodesic laminations on global
leaves can are mapped to the unit disk $\D$ by the identifications.
Thus we obtain a map from the local transverse set (obtained by
considering each local leaf in the chart as a point) to the space
of geodesic laminations on the unit disk $\D$. We require that
this map is continuous for the Hausdorff topology on the space of
geodesic laminations.

\vskip .2 cm

This definition certainly seems in the spirit of transverse
continuity of the hyperbolic metrics on $\S$. However, we
introduce below measured laminations on $\S$ in terms of the
continuity of measures. It turns out that the support of measured
laminations on $\S$ are not geodesic laminations as above, even
though the restriction to each leaf is a geodesic lamination.

\vskip .2 cm

\paragraph{\bf Definition 4.2} A (transversely continuous) {\it
measured lamination} $\mu$ on $\S$ is an assignment of a bounded
measured lamination to each leaf of $\S$ such that it is
continuous for the transverse variations with respect to Fr\'echet
topology on the space of measured laminations on the unit disk.

\vskip .1 cm

Given a local chart $D\times T$, where $D$ is a 2-disk and $T$ a
transverse Cantor set, the measured lamination $\mu$ on $\S$ gives
a map $\mu :T\to ML_{bdd}(\D )$ using the identifications of leaves
induced by the local chart. In the above definition, we require
that $\|\mu (t)-\mu (t')\|_{\nu}\to 0$ as $t'\to t$ for each $t\in
T$ and for each $0<\nu\leq 1$.

\vskip .2 cm

The definition of measured laminations on $\S$ does not specify
the support. To give an example of a measured lamination on $\S$,
fix a measured lamination $\sigma$ on a compact surface $S$ of
genus at least two. Since each leaf of $\S$ is a universal cover
of the compact surface $S$, we can lift $\sigma$ to a measured
lamination $\tilde{\sigma}$ on each leaf of $\S$. The lifts
$\tilde{\sigma}$ are locally constant for the transverse
variations in the local charts of a TLC hyperbolic metric coming
from the hyperbolic metric on $S$. Thus $\tilde{\sigma}$ defines a
geodesic lamination on $\S$.

\vskip .2 cm

As we mentioned above, there are measured laminations on $\S$
whose supports are not a geodesic laminations on $\S$. We give an
example of a such measured lamination.

\vskip .2 cm

\paragraph{\bf Example 4.3} We identify $\S$ with
$(\D\times\hat{G})/G$, for a Fuchsian group $G$ uniformizing a
compact hyperbolic surface $S=\D /G$. We define a measured
lamination $\tilde{\mu}$ on $\D\times\hat{G}$ which is invariant
under $G$. Let $G_i$ be a decreasing sequence of finite index
normal subgroups of $G=G_1$ such that $\cap_{i=1}^{\infty}G_i=\{
id\}$. We fix two simple closed curves $\gamma_1,\gamma_2$ in $S$
which intersect in one point and we fix two lifts
$\tilde{\gamma}_1,\tilde{\gamma}_2$ of $\gamma_1,\gamma_2$ in the
universal cover $\D$ such that
$|\tilde{\gamma}_1\cap\tilde{\gamma}_2|=1$. Denote by $C_1,C_2$
primitive hyperbolic translations in $G$ whose axes are
$\tilde{\gamma}_1,\tilde{\gamma}_2$ respectively. Let $C_1^{r_i}$
and $C_2^{t_i}$ be primitive elements in $G_i$. We further require
that the group $\tilde{G}_{i+1}$ generated by
$C_1^{r_i},C_2^{t_i}$ and $G_{i+1}$ is of index at least $3$ in
$G_i$.

\vskip .2 cm

Consider the cosets
$a_0^i\tilde{G}_{i+1},a_1^i\tilde{G}_{i+1},\ldots
,a_{k_i}^i\tilde{G}_{i+1}$, $a_0^i=id$, of $\tilde{G}_{i+1}$ in $G$.
Since $[G_i:\tilde{G}_{i+1}]\geq 3$, there are at least two cosets
different from $\tilde{G}_{i+1}$ which lie in $G_i$. Denote them by
$a_1^i\tilde{G}_{i+1},a_2^i\tilde{G}_{i+1}\subset G_i$. Then
$(a_1^iG_{i+1})\cdot (a_2^iG_{i+1})^{-1}$ are not of the form
$C_j^kG_{i+1}$, for $j=1$ or $j=2$ and for some $k\in\mathbb{Z}$. To
see this, first note that if $a_1^i(a_2^i)^{-1}$ is a power of $C_1$
or $C_2$ then it has to be a power of primitive elements $C_1^{r_i}$
or $C_2^{t_i}$ in $G_i$. (Otherwise $(a_1^iG_{i+1})\cdot
(a_2^iG_{i+1})^{-1}\notin G_i/G_{i+1}$ which is a contradiction.) On
the other hand, by our choice of $\tilde{G}_{i+1}$ and cosets
$a_1^i\tilde{G}_{i+1},a_2^i\tilde{G}_{i+1}$ we get that
$a_1^i(a_2^i)^{-1}$ is not a power of primitive elements $C_1^{r_i}$
or $C_2^{t_i}$.

\vskip .2 cm

Let
$\delta_{\tilde{\gamma}_i,t}$ denotes a unit mass measure on the
space of geodesics of $\D\times \hat{G}$ supported on the geodesic
$(\tilde{\gamma}_i,t)\subset \D\times\{ t\}$, for $t\in \hat{G}$.
We define
$$\tilde{\mu}':=\sum_i\Big{(} \sum_{t\in
a_1^i\hat{G}_{i+1}}m_i\delta_{\tilde{\gamma}_1,t}+\sum_{t\in
a_2^i\hat{G}_{i+1}}m_i\delta_{\tilde{\gamma}_2,t}\Big{)}
$$
where $m_i>0$, $m_i\to 0$ as $i\to\infty$ and
$\hat{G}_{i+1}<\hat{G}$ is the profinite completion of $G_{i+1}$.
The measured lamination $\tilde{\mu}'$ is varying continuously in
the transverse direction for the Fr\'echet topology on measured
laminations of the unit disk. The continuity is immediate at any
$t\in \hat{G}$ because $\tilde{\mu}'$ is locally constant.

\vskip .2 cm

We define
$$
\tilde{\mu}:=\sum_{A\in G}A^{*}(\tilde{\mu}').
$$
Then $\tilde{\mu}$ is invariant under the action of $G$. We first
show that the support of $\tilde{\mu}$ on each leaf $\D\times\{
t\}$, $t\in\hat{G}$, is a geodesic lamination. Let $\omega$ be a
fundamental polygon for the action of $G$ in $\D$. Then
$\omega\times\hat{G}$ is a fundamental set for the action of $G$
on $\D\times\hat{G}$. It is enough to show that the support of
$\tilde{\mu}$ has no self-intersections on $\omega\times\hat{G}$.
Note that the support of $\tilde{\mu}'$ has no self-intersections
because the cosets of $\hat{G}$ which contain copies of
$\tilde{\gamma}_1$ and $\tilde{\gamma}_2$ in the support are
chosen to be disjoint. The only possibility for the support of
$\tilde{\mu}$ to have a self-intersection is if a coset
$a_1^i\hat{G}_{i+1}$ is mapped onto the coset $a_2^i\hat{G}_{i+1}$
by $C_1^k$ for some $k\in\mathbb{Z}$, and similar for $C_2$. This
is impossible by our choice of cosets.

\vskip .2 cm

 We claim
that $\tilde{\mu}$ is continuous for the transverse variations in
the Fr\'echet topology on the space of measured laminations of the
unit disk $\D$. We first can assume that different lifts of
$\gamma_j$, $j=1,2$, in $\D$ do not belong to a single box of
geodesics $Q$ with $L(Q)=\log 2$ by appropriately choosing the group
$G$. If we show continuity in this case, the result follows because
the convergence of measured laminations in the Fr\'echet topology is
independent of the hyperbolic metric. We already concluded that
$\tilde{\mu}'$ is continuous for the transverse variations. By
taking the push-forward of $\tilde{\mu}'$ by $G$, we add some extra
support of $\tilde{\mu}$ intersecting fundamental set
$\omega\times\hat{G}$. The extra support is obtained by adding
$(\tilde{\gamma}_1,t)$ for $t\in C_1^ka_1^i\hat{G}_{i+1}$,
$k\in\mathbb{Z}$ and $i=1,2,\ldots$, and by adding
$(\tilde{\gamma}_2,t)$ for $t\in C_2^ka_2^i\hat{G}_{i+1}$,
$k\in\mathbb{Z}$ and $i=1,2,\ldots$. It is obvious that the
restriction of $\tilde{\mu}$ to the part which intersects
$\omega\times\hat{G}$ is continuous for the transverse variations at
any $t\in\hat{G}-\{ id\}$, similar to $\tilde{\mu}'$. The continuity
at $t=id$ follows because $m_i\to 0$. Our assumption that the orbit
of $\tilde{\gamma}_i$ does not contain two geodesic which lie in a
box $Q$ with $L(Q)=\log 2$ implies the continuity of $\tilde{\mu}$.

\vskip .2 cm

The measured lamination $\tilde{\mu}$ on $\D\times\hat{G}$
descends to a measured lamination $\mu$ on the universal
hyperbolic solenoid $\S$. The continuity of $\mu$ for the
transverse variations follows by the continuity of $\tilde{\mu}$.
It is clear that the support of $\mu$ is not a geodesic lamination
on $\S$ as in Definition 4.1 because it is not a closed set.
Moreover, the closure in $\S$ of the support of $\mu$ is not a
geodesic lamination because on the baseleaf it consists of the
full preimage of the two intersecting geodesics $\gamma_1$ and
$\gamma_2$ on the closed surface $S$. $\Box$

\section{Earthquake theorem}

We show that any two points in the Teichm\"uller space $T(\S )$ of
the universal hyperbolic solenoid $\S$ are connected by an
earthquake. We first need to recall certain facts from \cite{Sa3}
about arbitrary points in $T(\S )$.

\vskip .2 cm

A TLC solenoid $\S$ is homeomorphic to $(\D\times\hat{G})/G$, where
$\D$ is the unit disk. The space $\D\times\hat{G}$ is considered as
a universal cover of $\S$. Denote by $\pi :\D\times\hat{G}\to\S$ the
covering map. The action of $G$ is given by
$$
A(z,t):=(A(z),tA^{-1}),
$$
where $A$ can be considered as an element of $\hat{G}$.

\vskip .2 cm

A point in $T(\S )$ is given by a (differentiable) quasiconformal
map $f:\S\to X$, where $X$ is the universal hyperbolic solenoid with
an arbitrary hyperbolic metric (not necessarily TLC). We introduced
(see \cite{Sa3}) the universal (hyperbolic) cover to $X$ and the
covering group as follows. The action by $G$ does not introduce
identifications to the set $\{ 0\}\times\hat{G}$, $0\in\D$. Consider
a local chart $D\times T$ for $X$ which contains $f(\pi (\{
0\}\times\hat{G}))$ as a vertical set, where $T\equiv\hat{G}$ and
$D$ is a 2-disk with center at $0$. To fix the notation, we assume
that $f(\pi (\{ 0\}\times\hat{G}))$ corresponds to $\{ 0\}\times
T\subset D\times T$ in the local chart. We fix unit tangent vectors
at the points $f(\pi (\{ 0\}\times\hat{G}))$ corresponding under the
chart map to the unit tangent vectors at the points $\{ 0\}\times T$
along the positive axis in the chart $D\times T$.

\vskip .2 cm

The universal cover for $X$ is, by the definition, $\D\times T$ and
the covering map $\pi_X:\D\times T\to X$ is given by isometrically
mapping each $\D\times\{ t\}$ onto the leaf containing $f(\pi
(0,t))$ such that the origins are mapped onto the origins and the
unit tangent vectors along positive axes are mapped onto the unit
tangent vectors along the positive axes when considered in the chart
$D\times T$. The map $f:\S\to X$ lifts to a quasiconformal map
$$
\tilde{f}:\D\times\hat{G}\to\D\times T
$$
of the universal covers.

\vskip .2 cm

The action of $G$ is conjugated by $\tilde{f}$ to an action of a
group $G_X$ on $\D\times T$. An element $A$ of $G$ acts on
$\D\times\hat{G}$ by $A(z,t)=(A(z),tA^{-1})$, namely it acts on the
unit disk component by hyperbolic isometry and it shifts the leaf
$\D\times\{ t\}$ onto the leaf $\D\times\{ tA^{-1}\}$. In
particular, the action on the disk coordinate is independent of the
leaf (the second coordinate). We define a covering transformation
$A_f$ for $X$ by the formula
$$
(A_f\circ \tilde{f})(z,t)=(\tilde{f}\circ A)(z,t).
$$
An equivalent definition for $A_f$ is
$$
A_f(z,t):=(\pi_X)^{-1}(\pi_X(z,t),tA^{-1}),
$$
where $\pi_X(z,t)\in\S$ and $\pi_X(\cdot ,tA^{-1})$ stands for the
inverse of the covering map restricted to $\D\times\{ tA^{-1}\}$.
(Recall that the covering map when restricted to each leaf is an
isometry for the hyperbolic metric on leaves.) Thus $A_f$ is an
isometry on each leaf of $\D\times T$, but it varies with leaves.
The group of covering maps for the hyperbolic solenoid $X$ is
denoted by $G_X:=\tilde{f}G\tilde{f}^{-1}$. (If all $A_f$ in a
finite index subgroup of $G_X$ are constant in the transverse
direction then $X$ has a TLC hyperbolic metric.)

\vskip .3 cm

We show that earthquakes are transitive in $T(\S )$.

\vskip .2 cm

\paragraph{\bf Theorem 5.1} {\it A measured lamination $\mu$ on a solenoid $X$ with
an arbitrary hyperbolic metric gives an (leafwise) earthquake map
$E_{\mu}$ of $X$ into another solenoid $Y$ with hyperbolic metric
such that there exists a (differentiable) quasiconformal map
$f:X\to Y$ whose extension to the boundary of leaves coincides
with the extension of $E_{\mu}$. Any two points in the
Teichm\"uller space $T(\S )$ of the universal hyperbolic solenoid
$\S$ are connected by a unique earthquake along a measured
lamination on $\S$.}

\vskip .2 cm

\paragraph{\bf Proof} We first show that an earthquake map along a
measured lamination $\mu$ on the hyperbolic solenoid $X$ gives
another hyperbolic solenoid $Y$. We recall that $\mu$ is an
assignment of bounded measured laminations to the leaves of $X$
such that it varies continuously for the transverse variations.

\vskip .2 cm

Let $f:\S\to X$ be a differentiable quasiconformal map, where $\S
=(\D\times\hat{G})/G$. Recall the universal cover $\pi_X:\D\times
T\to X$ and lift $\mu$ to a measured lamination $\tilde{\mu}$ which
varies continuously for the transverse variations. In other words,
$\tilde{\mu}:T\to ML_{bdd}(\D )$ is continuous for the Fr\'echet
topology on the space $ML_{bdd}(\D )$ of bounded measured
laminations of $\D$ and it satisfies the invariance under the action
of $G_f$, i.e.
$$
(A_f(t))^{*}(\tilde{\mu}(t))=\tilde{\mu}(tA^{-1})
$$
for all $A_f\in G_f$, where $A_f(t)$ is the hyperbolic isometry of
$\D$ obtained by restricting $A_f$ to a map from $\D\times\{ t\}$
onto $\D\times\{ tA^{-1}\}$.

\vskip .2 cm

We consider a family of earthquakes $E_{\tilde{\mu}(t)}:\D\to\D$,
for $t\in T$, and we normalize them to fix $1$, $i$ and $-1$ on
the unit circle $S^1=\partial\D$. We first show that they induce a
family of quasisymmetric maps of $S^1\times T\equiv\partial
(\D)\times T$ onto itself which conjugate $G_f$ onto another group
of leafwise isometries. Let $h_t:=E_{\tilde{\mu}(t)}|_{S^1}$,
where $h_t:1,i,-1\mapsto 1,i,-1$. By the above invariance of
$\tilde{\mu}$ under the action of $A_f$ and by the fact that an
earthquake is determined, up to post-composition with a hyperbolic
isometry, by its measure \cite{Th}, we get
$$
h_{tA^{-1}}\circ A_f(t)=A_{*}(t)\circ h_t,
$$
where $A_{*}(t)$ is a hyperbolic isometry between $\D\times\{ t\}$
and $\D\times\{ tA^{-1}\}$ defined by the equation. To each
$A_f\in G_f$ we assign such $A_{*}$ which is a hyperbolic isometry
between leaves. The maps $A_{*}$, for all $A_f\in G_f$, form a
group $G_{*}$ isomorphic to $G_f$ under the conjugation by $h$.

\vskip .2 cm

We claim that $A_{*}(t)$ is continuous in $t$ for the standard topology
on the space of hyperbolic isometries of $\D$. To see this, it is enough to
show that the images of three fixed point on $S^1$ vary continuously in $t$. By
the above equation, we get
$$
A_{*}(t)=h_{tA^{-1}}\circ A_f(t)\circ h_t^{-1}.
$$
Since $\tilde{\mu}(t)$ varies continuously in $t$ and since
$h_t,h_{tA^{-1}}$ are properly normalized earthquakes, we conclude
that $h_t,h_{tA^{-1}}$ are continuous in $t$ for the topology of
pointwise convergence by \cite[Proposition 3.3]{Sa1}. By our assumption,
$A_f(t)$ is continuous in $t$. Thus $A_{*}(t)$ is continuous in $t$.

\vskip .2 cm

We need to show that the quotient $(\D\times T)/G_{*}$ is
quasiconformally equivalent to $X$, i.e. there exists a
homeomorphism $g:X\to (\D\times T)/G_{*}$ which is a
differentiable quasiconformal map on each leaf and which varies
continuously in the transverse direction (in the
$C^{\infty}$-topology). Let $g_t=ex(h_t:\D\to\D)$ be barycentric
extension of $h_t:S^1\to S^1$ (see \cite{DE} for the definition
and properties of barycentric extension). Recall that the family
$h_t$, $t\in T$, is continuous in the pointwise convergence
topology. Then the family $g_t$ of barycentric extensions is
continuous for the $C^{\infty}$-topology of $C^{\infty}$-maps
(over compact subsets of $\D$) by \cite{DE}.

\vskip .2 cm

We claim that $f_t$ is also continuous in the parameter $t$ for
the quasiconformal topology. (Note that the fact that the
earthquake measures converge in the Fr\'echet topology does not
imply that the extension of earthquakes to $S^1$ converge in the
quasisymmetric topology by the example in Section 3. At this point
we strongly use compactness of the solenoid $X$.) Recall that
$G_X$ has a compact fundamental set for the action on $\D\times T$
\cite{Sa4} (given by the image under $\tilde{f}$ of the
fundamental set $\omega\times\hat{G}$ for the action of $G$ on
$\D\times\hat{G}$). Thus the Beltrami coefficients of the family
$g_t$ are continuous in $t$ for the supremum norm over the
fundamental set of $X$. By the invariance of the quasisymmetric
family $h_t$ and by the conformal naturallity of barycentric
extension \cite{DE}, we obtain
$$
g_{tA^{-1}}\circ A_f(t)=A_{*}(t)\circ g_t,
$$
for $t\in T$. This invariance under $G_f$ and the continuity of
Beltrami coefficients of $g_t$ on the fundamental domain of $G_f$
implies that the Beltrami coefficients of $g_t$ are continuous in
$t$ for the essential supremum norm on the unit disk. Thus we
obtained a differentiable quasiconformal homeomorphism
$\tilde{g}:\D\times T\to\D\times T$, $\tilde{g}(\cdot
,t):=g_t(\cdot )$, which conjugates $G_f$ onto $G_{*}$.
Consequently it projects onto a quasiconformal homeomorphism
$g:X\to (\D\times T)/G_{*}$. Thus the earthquake $E_{\mu}$ defines
a new hyperbolic solenoid $Y:=(\D\times T)/G_{*}$ which is the
image of $X$. By its definition, the boundary values of each $g_t$
agree with $E_{\tilde{\mu}(t)}|_{S^1}$. This finishes the proof of
the first part of the Theorem.

\vskip .2 cm

It remains to show that any two points $[f:\S\to X],[g:\S\to Y]$
are connected by an earthquake along a measured lamination on $X$.
In other words, we need to find a measured lamination $\mu$ on $X$
such that $E_{\mu}$ maps $X$ onto $Y$ and that the extensions of
$E_{\mu}$ to the boundaries of leaves are equal to the extensions
of $g\circ f^{-1}$. We lift the maps $f$ and $g$ to the maps
$\tilde{f}:\D\times T\to \D\times T$ and $\tilde{g}:\D\times T\to
\D\times T$ of the universal covers of $X$ and $Y$. Let
$h_t:=\tilde{g}\circ\tilde{f}^{-1}|_{\D\times \{ t\}}$. Note that
the family of quasisymmetric maps $h_t$, $t\in T$, is continuous
in $t$ for the quasisymmetric topology. By Thurston's earthquake
theorem for the unit disk \cite{Th}, there exists a measured
lamination $\tilde{\mu}(t)$ such that
$E_{\tilde{\mu}(t)}|_{S^1}=h_t$. Moreover, $\tilde{\mu}(t)$ is a
bounded earthquake measure on $\D$. Since $h_t$ vary continuously,
we get that $\tilde{\mu}(t)$ vary continuously in $t$ for the
Fr\'echet topology by Proposition 3.2. The family $h_t$ satisfies
invariance properties with respect to $G_X$. Therefore, by the
uniqueness of earthquake measures \cite{Th}, the family of
corresponding earthquake measures $\tilde{\mu}(t)$ also satisfies
invariance properties. Thus it descend to the desired earthquake
measure $\mu$ on $X$. $\Box$

\vskip .2 cm

We recall that Proposition 3.2 states that if quasisymmetric maps
are close (in the quasisymmetric topology) then corresponding
earthquake measures are close (in the Fr\'echet topology). The
converse is false in general. However, we showed above that the
compactness of the universal hyperbolic solenoid $\S$ forces the
continuity of quasisymmetric maps on nearby leaves obtained by
earthquaking along transversely continuous measured laminations. The
proof extends along the same lines to show that if two measured
laminations on the solenoid $\S$ are close in the Fr\'echet topology
then the extension of earthquake maps to the boundary leaves are
close in the quasisymmetric topology. We obtained

\vskip .2 cm

\paragraph{\bf Corollary 5.2} {\it The earthquake map which assigns
to each bounded measured lamination on the universal hyperbolic
solenoid $\S$ the corresponding marked hyperbolic solenoid is a
homeomorphism between the space $ML(\S )$ of bounded measured
laminations and the Teichm\"uller space $T(\S )$.}

\section{Thurston's Boundary for $T(\S )$}

We recall the definition of the Liouville map $\L :T(\D )\to H(\D
)$ from the universal Teichm\"uller space $T(\D )$ to the space of
H\"older distributions $H(\D )$ of the unit disk $\D$. (Liouville
map first appears in \cite{Bo1} in the case of the Teichm\"uller
space of a compact Riemann surface and it is used in \cite{Sa3} to
introduce Thurston-type boundary to the universal Teichm\"uller
space $T(\D )$.)

\vskip .2 cm

The {\it universal Teichm\"uller space} $T(\D )$ is the set of all
quasisymmetric maps $h:S^1\to S^1$ which fix $1,i,-1$. The topology
on $T(\D )$ is defined by requiring that two quasisymmetric maps are
close if there exist their quasiconformal extensions to $\D$ whose
Beltrami coefficients are close in the essential supremum norm on
$\D$.

\vskip .2 cm

The {\it Liouville map } $\mathcal{L}:T(\D )\to H(\D )$ is defined
by taking the pull-back
$$
\mathcal{L}(h):=h_{*}(L)
$$
of the Liouville measure $L$ by the quasisymmetric maps $h\in T(\D
)$. The Liouville map $\mathcal{L}$ is a homemorphism of $T(\D )$
onto its image; the image $\mathcal{L}(T(\D ))\subset H(\D )$ is
closed and unbounded (see \cite{Sa3}). An {\it asymptotic ray} to
$L(T(\D ))$ is a path $t\Psi$, $t>0$ and $\Psi\in H(\D )$, such
that there exists a path $\alpha_t$, $t>0$, in $T(\D )$ with
$$
\frac{1}{t}\alpha_t\to\Psi,
$$
as $t\to\infty$. Each positive ray through the origin intersects the
image $L(T(\mathcal{G}))$ in at most one point. Therefore, under the
projection of the vector space $H(\D )$ to the unit sphere (in $H(\D
)$, for a fixed $\nu$-norm), the set $\L(T(\D ))$ is mapped
homeomorphically and its boundary corresponds to the asymptotic
rays. Thus, we consider asymptotic rays to $\mathcal{L}(T(\D ))$ as
a natural boundary to $T(\D )$. In \cite{Sa3}, we characterized the
boundary points of the universal Teichm\"uller space $T(\D )$ as all
asymptotic rays along bounded measured laminations. Namely,

\vskip .2 cm

\paragraph{\bf Theorem 6.1} \cite{Sa3} {\it The Liouville map
$\mathcal{L} :T(\D )\to H(\D )$ is a homeomorphism onto its image
and $\mathcal{L}(T(\D ))$ projects homeomorphically to the unit
sphere. The boundary of $T(\D )$ is identified by the above
embedding with the space of bounded projective measured
laminations $PML_{bdd}(\D )$. The (quasiconformal) mapping class
group $QMCG(\D )$ acts continuously on the closure $T(\D )\cup
PML_{bdd}(\D )$ of the universal Teichm\"uller space $T(\D )$.}

\vskip .4 cm

We introduce a Thurston-type boundary to the Teichm\"uller space
$T(\S )$ of the universal hyperbolic solenoid $\S$. The space of
geodesics on $\S =(\D\times\hat{G})/G$ is naturally identified with
the $G$-orbits of points in $(S^1\times S^1-diag)\times\hat{G}$
given by lifting a single geodesic on a leaf of $\S$ to the
universal cover $\D\times\hat{G}$. Since each leaf of $\S$ is
isometric to the hyperbolic plane, it supports the Liouville measure
on the space of its geodesics. Thus $\S$ has a leafwise Liouville
measure which lifts to a leafwise measure, called the {\it leafwise
Liouville measure} $L_{leaf}$, on the space of geodesics of the
universal cover $\D\times\hat{G}$.

\vskip .2 cm

A {\it (leafwise) H\"older distribution} $\Psi$ for the universal
hyperbolic solenoid $\S$ is a family of H\"older distributions
$\Psi_t\in H(\D )$, for $t\in\hat{G}$, which are invariant under the
action of $G$, i.e.
$$
\Psi_{tA^{-1}}(\varphi\circ A^{-1})=\Psi_t(\varphi),
$$
where $\varphi$ is a H\"older continuous function with compact
support on the space of geodesics $\mathcal{G}(\D )$, and which vary
continuously in $t$ for the Fr\'echet topology, i.e.
$$
\|\Psi_t-\Psi_{t_1}\|_{\nu}\to 0,
$$
as $t_1\to t$ for each $t\in\hat{G}$ and for each $0<\nu\leq 1$.

\vskip .2 cm

The $\nu$-{\it norm} of a leafwise H\"older distribution $\Psi$ is
given by
$$
\|\Psi\|_{\nu}:=\sup_{t\in\hat{G},\varphi\in test(\nu
)}|\Psi_t(\varphi )|,
$$
for $0<\nu\leq 1$, where $test(\nu )$ is the set of $\nu$-test
functions on $\G (\D )$. If $\|\Psi\|_{\nu}<\infty$ for all
$0<\nu\leq 1$ then $\Psi$ is called {\it bounded leafwise H\"older
distribution}. The space of all (bounded) leafwise H\"older
distributions for the universal hyperbolic solenoid $\S$ is denoted
by $H(\S )$.

\vskip .2 cm

Let $[f:\S\to X]\in T(\S )$ be an arbitrary point and
denote by $\tilde{f}:\D\times\hat{G}\to\D\times T$ the lift of
$f$ to the universal cover. Let $h:S^1\times \hat{G}\to S^1\times T$
be the leafwise quasisymmetric extension of $\tilde{f}$ to the boundary of leaves. We define
the {\it Liouville map} $\mathcal{L}_{\S}:T(\S )\to H(\S )$ for the universal hyperbolic solenoid $\S$ by
$$
\mathcal{L}_{\S}([f])=h_{*}(L_{leaf}).
$$
(Note that bounded measures on $\mathcal{G}(\D )$ are in $H(\D )$
and that a pull-back by a quasisymmetric of the Liouville measure
is bounded \cite{Sa3}. Thus the image of the Liouville map is in
$H(\D )$ and the leafwise statement for the universal hyperbolic
solenoid immediately follows.)

\vskip .2 cm

We show that the Liouville map is an embedding and that the
natural boundary (i.e. the set of asymptotic rays) is homeomorphic
to the space of projective measured laminations on $\S$.

\vskip .2 cm

\paragraph{\bf Theorem 6.2} {\it The Liouville map
$\mathcal{L}_{\S}:T(\S )\to H(\S )$ is a homeomorphisms onto
its image. The set of asymptotic rays to $\mathcal{L}_{\S}(T(\S ))$ is
homeomorphic to the space of projective measured laminations on $\S$. The
baseleaf preserving mapping class group $MCG_{BLP}(\S )$ acts continuously
on the closure $T(\S )\cup PML(\S )$ of the Teichm\"uller space $T(\S )$
of the universal hyperbolic solenoid $\S$.}

\vskip .2 cm

\paragraph{\bf Proof}
The Liouville map $\L_{\S}$ is assigning to any $[f]\in T(\S )$ the
pull-backs of the leafwise Liouville measures on
$\mathcal{G}(\D\times\hat{G})$ by the extensions $h$ to $S^1\times
\hat{G}$ of the lift $\tilde{f}:\D\times\hat{G}\to\D\times\hat{G}$.
The continuity of $\L :T(\D )\to H(\D )$ implies that $(h_t)_{*}(L)$
is continuous in $t$ for the Fr\'echet topology. Thus $\L$ maps
$T(\S )$ into $H(\S )$. Recall that $T(\S )$ embeds in the universal
Teichm\"uller space $T(\D )$ by restricting the map $f:\S\to X$ to
the baseleaf of $\S$ \cite{NS}. Denote by $T_{restr.}(\D )$ the
image of the embedding. Also, since the baseleaf is dense in $\S$,
the restriction to the beaseleaf of the pull-back of the leafwise
Liouville measure completely determines the measure. Therefore, the
restriction of the Liouville map $\mathcal{L}_{\S}$ to the
baseleaf($\equiv\D$) $\mathcal{L}:T_{restr.}(\D ) \to H(\D )$
completely determines the map. Since $\L :T(\D )\to H(\D )$ is a
homeomorphism onto its image, it follows that
$\mathcal{L}:T_{restr.}(\D ) \to H(\D )$ is also a homeomorphism
onto its image. Therefore, $\mathcal{L}_{\S}:T(\S )\to H(\S )$ is a
homeomorphisms onto its image.

\vskip .2 cm

Let $s\beta$, $s>0$, be an asymptotic ray to $\L (T(\S ))$ in
$H(\S )$, i.e. there exists a path $s\mapsto\alpha_s$ in $\L (T(\S
))$ such that the lifts $\tilde{\alpha}_s,\tilde{\beta}$ to the
universal cover $\D\times\hat{G}$ satisfy
$\|\frac{1}{s}\tilde{\alpha}_s-\tilde{\beta}\|_{\nu}\to 0$ as
$s\to\infty$, for all $0<\nu\leq 1$. This implies that
$\|\tilde{\alpha}_s|_{\D\times\{ t\}}-\tilde{\beta}|_{\D\times\{
t\}}\|_{\nu}\to 0$ for all $t\in\hat{G}$ as $s\to\infty$. By
\cite{Sa3}, each $\tilde{\beta}|_{\D\times\{ t\}}$, $t\in\hat{G}$,
is a measured lamination on $\D\times\{ t\}$. Therefore, $\beta\in
H(\S )$ is a leafwise measured lamination, and since it belongs to
$H(\S )$, it is continuous for the transverse variations. Thus
$\beta$ is a measured lamination on $\S$. The earthquake theorem
(Theorem 5.1) for the universal hyperbolic solenoid $\S$ shows
that an earthquake path $E_{s\mu}:\S \to X_s$, $s>0$, is in $T(\S
)$. Denote by $\tilde{\mu}$ the lift of $\mu$ to the universal
cover $\D\times\hat{G}$ of $\S$. Then $h_s:=E_{s\tilde{\mu}}$ is a
leafwise earthquake map. By \cite[Theorem 2]{Sa3}, we have
$$
\frac{1}{s}(h_s|_{\D\times\{ t\}})_{*}(L)\to\tilde{\mu}|_{\D\times\{ t\}},
$$
in the $\nu$-norm, $0<\nu\leq 1$, as $s\to\infty$ for each
$t\in\hat{G}$. Moreover, \cite[Lemma 4.4]{Sa3} shows that the
above convergence in the $\nu$-norm, $0<\nu\leq 1$, is uniform
independent of the leaf $\D\times\{ t\}$. Thus
$$
(h_s)_{*}(L_{leaf})\to \mu,
$$
as $s\to\infty$ in the Fr\'echet topology on $H(\S )$. Thus the
natural boundary to $T(\S )$ (i.e. the space of asymptotic rays to
$\L (T(\S ))$ in $H(\S )$) is homeomorphic to the space $PML(\S )$
of projective measured laminations on $\S$. The continuity of the
action of the baseleaf preserving mapping class group on the
closure $T(\S )\cup PML(\S )$ is immediate from \cite{Sa3}. $\Box$

\vskip .4 cm

Note that leafwise measured laminations on $\S$ are given in terms
of continuity for the transverse variations. We show that each
leafwise measured lamination on $\S$ is approximated by
transversely locally constant (TLC) measured lamination (i.e.
laminations obtained by lifting laminations on compact surfaces to
$\S$) which is parallel to the statement that each hyperbolic
metric on $\S$ is approximated by TLC hyperbolic metrics.

\vskip .2 cm

\paragraph{\bf Theorem 6.3} {\it The subset of all measured
lamination on the universal hyperbolic solenoid $\S$ which are
locally transversely constant (TLC) is dense in the space of all
measured laminations $ML(\S )$ on $\S$ for the Fr\'echet
topology.}

\vskip .2 cm

\paragraph{\bf Proof}
Recall that $T_{restr.}(\D )$ is obtained by taking the closure of
the union of all quasisymmetric maps which conjugate finite index
subgroups of $G$ onto other Fuchsian groups \cite{NS}, i.e. the
closure of the union $\cup_{[G:K]<\infty}\tilde{T}(\D /K)$ of the
lifts $\tilde{T}(\D /K)$ to $T(\D )$ of all Teichm\"uller spaces
$T(\D /K)$ of finite degree unbranched covers $\D /K$ of $\D /G$. By
the characterization of the image of $\L$ from \cite{Bo1}, the image
$\mathcal{L}(\cup_{[G:K]<\infty}\tilde{T}(\D /K))$ of the union
consists of all bounded H\"older distributions $\alpha$ which are
positive measures invariant under finite index subgroups of $G$ and
which satisfy $e^{-\alpha ([a,b]\times [c,d])}+e^{-\alpha
([b,c]\times [d,a])}=1$ for all $a,b,c,d\in S^1$ given in
counterclockwise order.

 \vskip .2 cm

Since $\L$ is a homeomorphism onto its image and $\L (T(\D ))$ is
closed in $H(\D )$ \cite{Sa3}, it follows that $\L (T_{restr.}(\D
))$ equals to the closure (in the Fr\'echet topology) of $\L
(\cup_{[G:K]<\infty}\tilde{T}(\D /K))$. By \cite{Bo1} or \cite{Sa3},
the asymptotic rays to $\L (\cup_{[G:K]<\infty}\tilde{T}(\D /K))$
 are containing all the
projective measured laminations which are invariant under all
finite index subgroups $K$ of $G$. The Liouville map $\L$ composed
with the projection $pr:H(\S )\to S_{\nu}$ to the unit sphere
$S_{\nu}$ in $H(\S )$ (for a fixed $\nu$-norm) is a homeomorphism
(see \cite{Sa3}). This implies that the points in the closure of
$pr\circ\L (T_{restr.}(\D ))$ which are not in $pr\circ\L
(T_{restr.}(\D ))$ are projectivized asymptotic rays to $\L
(T_{restr.}(\D )$. Any such point in the closure of $pr\circ\L
(T_{restr.}(\D ))$ is approximated by projective measures in
$pr\circ\L (T_{restr.}(\D ))$ that are invariant under finite
index subgroups of $G$. Thus the closure of all invariant (under
finite index subgroups of $G$) projective measured laminations
contains all asymptotic rays to $\L (T_{restr.})$. The measured
laminations invariant under finite index subgroups of $G$ lift to
locally transversely constant measured laminations on $\S$. Thus
the limits of locally transversely constant measured laminations
on $\S$ give all (transversely continuous) measured laminations on
$\S$. $\Box$

\vskip .2 cm

\paragraph{\bf Remark 6.4} The set of asymptotic rays to $\L
(T_{restr.}(\D ))$ in $H(\D )$ is equal to the restriction to the
baseleaf of the set of asymptotic rays to $\L (T(\S ))$ in $H(\S )$.
This is a consequence of the proof of Theorem 6.3, since each
asymptotic ray for $\L (T_{restr.}(\D ))$ is the limit in the
Fr\'echet topology of the asymptotic rays invariant under finite
index subgroups of $G$.

\section{The punctured solenoid}

We sketch an extension of our results to the punctured solenoid
$\S_p$ defined in \cite{PS}. We first recall the definition of $\S_p$.

\vskip .2 cm

 Let $H$ be the subgroup of $PSL_2(\mathbb{Z} )$ such that $\D /H$
 is a once punctured torus. Denote by $\hat{H}$ the profinite
 completion of $H$. Then we define (see \cite{PS}) the punctured solenoid by
 $$
 \S_p:=(\D\times\hat{H})/H,
 $$
 where the action of $H$ is given by $A(z,t):=(A(z),tA^{-1})$, for
 $A\in H$.

\vskip .2 cm

A leafwise measured lamination on $\S_p$ is an assignment of a
bounded measured lamination to each leaf of $\S_p$ which varies
continuously in the transverse direction. The support of a
leafwise measured lamination on $\S_p$ is a leafwise geodesic
lamination which is not necessarily continuous for the transverse
variations. We say that the support of a leafwise measured
lamination on $\S_p$ is compact if, when restricted to each leaf
of $\S_p$, the support geodesic lamination is a pre-compact subset
of $\S_p$. The earthquake theorem holds for $\S_p$ when we use the
space of measured laminations with compact support $ML_0(\S_p)$.

\vskip .2 cm

If a simple geodesic on a punctured surface enters a definite
neighborhood of a puncture, then it has an endpoint at the
puncture. We recall a standard proof of this fact in the upper
half-plane model of the hyperbolic plane. Let $A:z\mapsto z+1$ be
the parabolic element corresponding to the puncture. If a lift of
the geodesic entering a neighborhood of a puncture on the surface
is a Euclidean half-circle with radius greater than $1/2$, then
the translate of the lift of the geodesic by $A$ intersects the
lift of the geodesic. Thus the geodesic is not simple.
Contradiction. Therefore, there exists a definite neighborhood of
a puncture where no simple geodesic enters unless it ends at the
puncture.

\vskip .2 cm

Assume that a measured lamination $\mu$ on $\S_p$ does not have a
compact support. Then there exists a leaf $l$ of $\S_p$ such that
the restriction $\mu_l$ of $\mu$ contains a geodesic $g\in l$ with
an endpoint at the puncture. By the continuity for the transverse
variations, the support of $\mu_l$ contains the translates of $g$ in
$l$ by parabolic elements $A^{kn}\in H$, for a fixed $k>0$ and for
all $n\in\mathbb{N}$, where $A$ fixes the endpoint of $g$. Then $g$
has to be isolated in $l$ because otherwise the translates under
$A^{kn}$ of the geodesics in the support of $\mu_l$ converging to
$g$ would intersect the geodesics converging to $g$ similar to the
punctured surface case. Thus $\mu_l$ would not be a measured
lamination. The remaining possibility is that $g$ is isolated with
atomic measure. Then the translates $A^{kn}(g)$ also have atomic
measure approximately equal to the atomic measure of $g$ for $n$
large and they share the same endpoint. Then the measured lamination
$\mu|_{l}$ is not bounded which is a contradiction with the
definition of a measured lamination on $\S_p$. Therefore, a measured
lamination $\mu$ on $\S_p$ always has a compact support.

\vskip .2 cm

As in the proof of Theorem 5.1 it follows that the extension to
the boundary of each quasiconformal map
$f:\S_p\to X_p$ can be achieved by a leafwise earthquakes whose measures vary
continuously in the transverse direction for the Fr\'echet topology on earthquake
measures. On the other hand, a measured lamination $\mu$ on $\S_p$ with compact support
gives a quasiconformal map $f:\S_p\to X_p$ whose extensions to the boundaries of the
leaves agrees with the earthquake $E_{\mu}$. Thus the earthquake theorem is true for the
punctured solenoid for the space of measured laminations $ML_0(\S_p )$ with the compact support
on $\S_p$.

\vskip .4 cm

We explain how to extend Thurston's boundary to the Teichm\"uller
space of the punctured solenoid $\S_p$. We show that the boundary
consists of $PML_0(\S_p)$. By the extension of the earthquake
theorem and by the proof of Theorem 6.2, it is only necessary to
show that the asymptotic rays to $\L (T(\S_p ))\subset H(\S_p)$
are of the form $tW$, $t>0$, for some $W\in ML_0(\S_p)$. Let
$\alpha_t\in\L (T(\S_p))$ be such that $\frac{1}{t}\alpha_t\to W$
as $t\to\infty$. Then the restriction of $W$ to each leaf of
$\S_p$ is a bounded measured lamination by the results in
\cite{Sa3}. Since $W$ is continuous for the transverse variations,
it follows that $W$ is a measured lamination on $\S_p$. By the
above, $W\in ML_0(\S_p)$.

\end{document}